\documentclass[12pt,a4paper,reqno]{amsart} %Fr o m detta
\usepackage{amssymb}
\usepackage{latexsym}
\usepackage{amsmath}
\usepackage{mathrsfs}
\usepackage[X2,T1]{fontenc}
\usepackage[applemac]{inputenc}
\usepackage{cite}
\usepackage{calc}                   %Detta beh?vs f?r att
\newcommand{\scal}[2]{\langle #1,#2\rangle}
\newcommand{\rr}[1]{\mathbf R^{#1}}
\newcommand{\cc}[1]{\mathbf C^{#1}}
\newcommand{\nm}[2]{\Vert #1\Vert _{#2}}

\newcommand{\sets}[2]{\{ \, #1\, ;\, #2\, \} }

\newcommand{\ep}{\varepsilon}

\newcommand{\cdo}{\, \cdot \, }

\newcommand{\eabs}[1]{\langle #1\rangle}     %%%%%   for <x>

\newcommand{\vrum}{\vspace{0.1cm}}

\newcommand{\GS}{\operatorname{GS}}

\newcommand{\nn}[1]{{\mathbf N}^{#1}}
\newcommand{\maclA}{\mathcal A}
\newcommand{\maclS}{\mathcal S}
\newcommand{\mascF}{\mathscr F}
\newcommand{\mascS}{\mathscr S}

\def\sS{{\mathscr S}}
\def\cS{{\mathcal S}}

\def\R{\mathbf{R}}

\numberwithin{equation}{section}          %Detta g?r att man f?r
\newtheorem{thm}{Theorem}
\numberwithin{thm}{section}

\newcommand{\rubrik}{}
\newtheorem{prop}[thm]{Proposition}

\newtheorem{lemma}[thm]{Lemma}

\theoremstyle{definition}

\theoremstyle{remark}
\newtheorem{rem}[thm]{Remark}

\author{Marco Cappiello}

\address{Dipartimento di Matematica ``G. Peano", Universit\`a di Torino,
Via Carlo Alberto 10, Torino, Italy}

\email{marco.cappiello@unito.it}

\author{Luigi Rodino}

\address{Dipartimento di Matematica ``G. Peano", Universit\`a di Torino,
Via Carlo Alberto 10, Torino, Italy}

\email{luigi.rodino@unito.it}

\author{Joachim Toft}

\address{Department of Mathematics,
Linn{\ae}us University, V{\"a}xj{\"o}, Sweden}

\email{joachim.toft@lnu.se}

\title{Radial symmetric elements and the Bargmann transform}



\begin{document}

\begin{abstract}
We prove that a function or distribution on $\rr d$
is radial symmetric, if and only if its Bargmann transform
is a composition by an entire function on $\mathbf C$
and the canonical quadratic function from $\cc d$ to
$\mathbf C$. 
\end{abstract}

\maketitle

\par

%%%%%%%%%%%%%%%%%%%%%%%
\section{Introduction}\label{sec0}
%%%%%%%%%%%%%%%%%%%%%%%

\par

A function or distribution on $\rr d$ is called radial symmetric if it is
invariant under pullbacks of unitary transformations on $\rr d$. In
the paper we prove that any function or distribution $f$ on $\rr d$ is
radial symmetric, if and only if its Bargmann transform $\mathfrak Vf$
satisfy
$$
\mathfrak Vf (z) = F_0(\scal zz), \quad z\in \cc d,
$$
for some entire function $F_0$ on $\mathbf C$. We also prove that
$f$ is radial symmetric is equivalent to that $f$ is orthogonal in $L^2(\R^d)$
to every Hermite function $h_\alpha$, as long as at least one of
$\alpha _j$ is odd, and that
$$
\frac {\alpha !}{\sqrt {(2\alpha )!}}(f,h_{2\alpha})
$$
only depends on $|\alpha|$.

\par

We perform these investigations in the framework of Schwartz or
Gelfand-Shilov functions, and corresponding distribution spaces.
To this end we devote the preliminary part of the paper to results 
on the Bargmann transform on these spaces.

\par

We also use our results to show that there is a natural way
to assign to any radial symmetric function or distribution $f$ 
on $\rr d$,  a distribution on $\mathbf R$ which obeys similar
estimates as $f$.

\par

%%%%%%%%%%%%%%%%%%%%%%%
\section{Preliminaries}\label{sec1}
%%%%%%%%%%%%%%%%%%%%%%%

\par

In this section we recall some basic properties on the
Bargmann transform. We shall often formulate these results in the
framework of the Gelfand-Shilov space $\cS _{1/2}(\rr d)$ and
its dual $\cS _{1/2}'(\rr d)$ (see e.{\,}g. \cite{GS}). The reader
who is not interested in this general situation may replace
$\cS _{1/2}(\rr d)$ and $\cS _{1/2}'(\rr d)$ by $\sS (\rr d)$
and $\sS '(\rr d)$ respectively. Here $\sS (\rr d)$ is the set
of Schwartz functions on $\rr d$, and $\sS '(\rr d)$ is the
set of tempered distributions on $\rr d$, see for example
\cite{Ho1}.

\par

\subsection{Gelfand-Shilov spaces}\label{subsec1.2}
We start by recalling some facts about Gelfand-Shilov spaces.
Let $0<h,s,t\in \mathbf R$ be fixed. Then $\mathcal S_{s,h}(\rr d)$
consists of all $f\in C^\infty (\rr d)$ such that
\begin{equation*}%\label{gfseminorm}
\nm f{\mathcal S_{t,h}^s}\equiv \sup \frac {|x^\beta \partial ^\alpha
f(x)|}{h^{|\alpha | + |\beta |}\alpha !^s\, \beta !^t}
\end{equation*}
is finite. Here the supremum should be taken over all $\alpha ,\beta \in
\mathbf N^d$ and $x\in \rr d$.

\par

Obviously $\mathcal S_{s,h}^t$ is a Banach space which increases with $h$, $s$ and $t$ and 
$\mathcal S_{s,h}^t\hookrightarrow \mathscr S$. Here and
in what follows we use the notation $A\hookrightarrow B$ when the topological
spaces $A$ and $B$ satisfy $A\subseteq B$ with continuous embeddings.
Furthermore, if $s+t\ge 1$ and $(s,t)\neq (1/2,1/2)$, or $s =t=1/2$ and $h$
is sufficiently large, then $\mathcal
S_{t,h}^s$ contains all finite linear combinations of Hermite functions.
Since such linear combinations are dense in $\mathscr S$, it follows
that the dual $(\mathcal S_{t,h}^s)'(\rr d)$ of $\mathcal S_{t,h}^s(\rr d)$ is
a Banach space which contains $\mathscr S'(\rr d)$.

\par

The \emph{Gelfand-Shilov spaces} $\mathcal S_{t}^s(\rr d)$ and
$\Sigma _{t}^s(\rr d)$ are defined as the inductive and projective 
limits respectively of $\mathcal S_{t,h}^s(\rr d)$. This implies that
\begin{equation}\label{GSspacecond1}
\mathcal S_t^{s}(\rr d) = \bigcup _{h>0}\mathcal S_{t,h}^s(\rr d)
\quad \text{and}\quad \Sigma _t^{s}(\rr d) =\bigcap _{h>0}
\mathcal S_{t,h}^s(\rr d),
\end{equation}
and that the topology for $\mathcal S_t^{s}(\rr d)$ is the strongest
possible one such that the inclusion map from $\mathcal S_{t,h}^s
(\rr d)$ to $\mathcal S_t^{s}(\rr d)$ is continuous, for every choice 
of $h>0$. The space $\Sigma _t^s(\rr d)$ is a Fr{\'e}chet space
with seminorms $\nm \cdo{\mathcal S_{t,h}^s}$, $h>0$. Moreover,
$\Sigma _t^s(\rr d)\neq \{ 0\}$, if and only if $s+t\ge 1$ and
$(s,t)\neq (1/2,1/2)$, and $\maclS _t^s(\rr d)\neq \{ 0\}$, if and only
if $s+t\ge 1$. From now on we
assume that the Gelfand-Shilov parameter pair $(s,t)$ are
\emph{admissible}, or \emph{$\GS$-admissible},
that is, $s+t\ge 1$ and $(s,t)\neq (1/2,1/2)$ when considering
$\Sigma _t^s(\rr d)$, and $s+t\ge 1$ when considering
$\maclS _t^s(\rr d)$

\medspace

The \emph{Gelfand-Shilov distribution spaces} $(\mathcal S_t^{s})'(\rr d)$
and $(\Sigma _t^s)'(\rr d)$ are the projective and inductive limit
respectively of $(\mathcal S_{t,h}^s)'(\rr d)$.  This means that
\begin{equation}\tag*{(\ref{GSspacecond1})$'$}
(\mathcal S_t^s)'(\rr d) = \bigcap _{h>0}(\mathcal S_{t,h}^s)'(\rr d)\quad
\text{and}\quad (\Sigma _t^s)'(\rr d) =\bigcup _{h>0}(\mathcal S_{t,h}^s)'(\rr d).
\end{equation}
We remark that in \cite{GS, Pil} it is proved that $(\mathcal S_t^s)'(\rr d)$
is the dual of $\mathcal S_t^s(\rr d)$, and $(\Sigma _t^s)'(\rr d)$
is the dual of $\Sigma _t^s(\rr d)$ (also in topological sense). For conveniency we
set
$$
\maclS _s=\maclS _s^s,\quad \maclS _s'=(\maclS _s^s)',\quad
\Sigma _s=\Sigma _s^s,\quad \Sigma _s'=(\Sigma _s^s)'.
$$

\par

For every admissible $s,t>0$ and $\ep >0$ we have
\begin{equation}\label{GSembeddings}
\begin{alignedat}{2}
\Sigma _t^s (\rr d) &\hookrightarrow &
\maclS _t^s(\rr d) &\hookrightarrow  \Sigma _{t+\ep}^{s+\ep}(\rr d)
\\[1ex]
\quad \text{and}\quad
(\Sigma _{t+\ep}^{s+\ep})' (\rr d) &\hookrightarrow & (\maclS _t^s)'(\rr d)
&\hookrightarrow  (\Sigma _t^s)'(\rr d).
\end{alignedat}
\end{equation}

\par

The Gelfand-Shilov spaces possess several convenient mapping
properties, and in the case $s=t$ they are invariant under
several basic transformations. For example they are invariant under
translations, dilations, tensor products
and under (partial) Fourier transformations.

\par

From now on we let $\mathscr F$ be the Fourier transform which
takes the form
$$
(\mathscr Ff)(\xi )= \widehat f(\xi ) \equiv (2\pi )^{-d/2}\int _{\rr
{d}} f(x)e^{-i\scal  x\xi }\, dx
$$
when $f\in L^1(\rr d)$. Here $\scal \cdo \cdo$ denotes the usual
scalar product on $\rr d$. The map $\mathscr F$ extends 
uniquely to homeomorphisms on $\mathscr S'(\rr d)$,
from $(\mathcal S_t^s)'(\rr d)$ to $(\mathcal S_s^t)'(\rr d)$ and
from $(\Sigma _t^s)'(\rr d)$ to $(\Sigma _s^t)'(\rr d)$. Furthermore,
$\mascF$ restricts to
homeomorphisms on $\mathscr S(\rr d)$, from
$\mathcal S_t^s(\rr d)$ to $\mathcal S_s^t(\rr d)$ and
from $\Sigma _t^s(\rr d)$ to $\Sigma _s^t(\rr d)$,
and to a unitary operator on $L^2(\rr d)$.

\par

It follows from the following lemma that elements in Gelfand-Shilov
spaces can be characterized by estimates of the form
\begin{equation}\label{GSexpcond}
|f(x)|\lesssim e^{-\ep |x|^{1/t}}\quad \text{and}\quad |\widehat f (\xi )|
\lesssim e^{-\ep |\xi |^{1/s}} .
\end{equation}
The proof is omitted, since the result can be found in e.{\,}g.
\cite{ChuChuKim, NiRo}. Here and in the sequel, $A\lesssim B$
means that $A\le cB$ for a suitable constant $c>0$.

\par

\begin{lemma}\label{GSFourierest}
Let $s,t>0$ and $f\in \mathcal S'_{1/2}(\rr d)$. Then the following is true:
\begin{enumerate}
\item[(1)] if $s+t\ge 1$, then $f\in \mathcal S_t^s(\rr d)$, if and
only if \eqref{GSexpcond} holds for some $\ep >0$;

\vrum

\item[(2)]  if $s+t\ge 1$ and $(s,t)\neq (1/2,1/2)$, then
$f\in \Sigma _t^s(\rr d)$, if and only
if \eqref{GSexpcond} holds for any $\ep >0$.
\end{enumerate}
\end{lemma}

\par

Gelfand-Shilov spaces and their distribution spaces can also,
in some sense more convenient ways, be
characterized by means of estimates of short-time Fourier transforms,
(see e.{\,}g. \cite{SiTo, Teof}).
We recall here the details and start by recalling the definition of
the short-time Fourier transform.

\par

Let $\phi \in \maclS _s '(\rr d)$ be fixed. Then the \emph{short-time
Fourier transform} $V_\phi f$ of $f\in \maclS _s '
(\rr d)$ with respect to the \emph{window function} $\phi$ is
the Gelfand-Shilov distribution on $\rr {2d}$, defined by
$$
V_\phi f(x,\xi ) \equiv  (\mascF _2 (U(f\otimes \phi )))(x,\xi ) =
\mascF (f \, \overline {\phi (\cdo -x)})(\xi
),
$$
where $(UF)(x,y)=F(y,y-x)$. If $f ,\phi \in \maclS _s (\rr d)$, then it follows that
$$
V_\phi f(x,\xi ) = (2\pi )^{-d/2}\int f(y)\overline {\phi
(y-x)}e^{-i\scal y\xi}\, dy .
$$

\par

The next two results show that both spaces of Gelfand-Shilov
functions and Gelfand-Shilov distributions can be completely
identified with growth and decay properties of the short-time
Fourier transforms for the involved functions and distributions.
The conditions are of the forms
\begin{align}
|V_\phi f(x,\xi )| &\lesssim  e^{-\ep (|x|^{1/t}+|\xi |^{1/s})},
\label{stftexpest2}
\\[1ex]
|(\mathscr F(V_\phi f))(\xi ,x)| &\lesssim  e^{-\ep (|x|^{1/t}+|\xi |^{1/s})}
\label{stftexpest3}
\intertext{and}
|V_\phi f(x,\xi )| &\lesssim  e^{\ep (|x|^{1/t}+|\xi |^{1/s})}.
\tag*{(\ref{stftexpest2})$'$}
\end{align}

\par

\begin{prop}\label{stftGelfand2}
Let $(s,t)$ and $(s_0,t_0)$ be $\GS$-admissible which satisfy
$s_0\le s$ and $t_0\le t$, and let
$\phi \in \mathcal S_{t_0}^{s_0}(\rr d)\setminus 0$ and
$f\in (\mathcal S_{t_0}^{s_0})'(\rr d)$.
Then the following is true:
\begin{enumerate}
\item $f\in  \mathcal S_{t}^s(\rr d)$, if and only if
\eqref{stftexpest2} holds for some $\ep > 0$;
%Furthermore, if $f\in  \mathcal S_{t}^s(\rr d)$, then
%\eqref{stftexpest3} holds for some constants $\ep > 0$;

\vrum

\item if in addition $\phi \in \Sigma _{t}^{s}(\rr d)$, then
$f\in  \Sigma _{t}^s(\rr d)$, if and only if \eqref{stftexpest2}
holds for every $\ep > 0$.
\end{enumerate}
\end{prop}

\par

A proof of Theorem \ref{stftGelfand2} can be found in e.{\,}g. \cite{GZ}
(cf. \cite[Theorem 2.7]{GZ}). The corresponding result for Gelfand-Shilov
distributions is the following, which is essentially a restatement of
\cite[Theorem 2.5]{To11}.

\par

\begin{prop}\label{stftGelfand2dist}
Let $(s,t)$ and $(s_0,t_0)$ be $\GS$-admissible which satisfy
$1/2<s_0\le s$ and $1/2<t_0\le t$, and let $\phi \in \Sigma _{t}^s(\rr d)
\setminus 0$, and let $f\in (\mathcal S_{t_0}^{s_0})'(\rr d)$. Then
the following is true:
\begin{enumerate}
\item $f\in  (\mathcal S_{t}^s)'(\rr d)$, if and only if
\eqref{stftexpest2}$'$ holds for every
$\ep > 0$;

\vrum

\item $f\in  (\Sigma _{t}^s)'(\rr d)$, if and only if \eqref{stftexpest2}$'$
holds for some $\ep > 0$.
\end{enumerate}
\end{prop}

\par

We note that in (2) in \cite[Theorem 2.5]{To11} it should stay $(\Sigma _{t}^s)'(\rr d)$
instead of $\Sigma _{t}^s(\rr d)$.

\par

\begin{rem}\label{SchwFunctionSTFT}
The short-time Fourier transform can also be used to identify the elements
in $\mascS (\rr d)$ and in $\mascS '(\rr d)$. In fact, if
$\phi \in \mascS (\rr d)\setminus 0$ and $f\in (\maclS _{1/2})'(\rr d)$,
then the following is true:
\begin{enumerate}
\item $f\in \mascS (\rr d)$, if and only if for every $N\ge 0$, it holds
$$
|V_\phi f(x,\xi )| \lesssim \eabs {x,\xi }^{-N} \text ;
$$

\vrum

\item $f\in \mascS '(\rr d)$, if and only if for some $N\ge 0$, it holds
$$
|V_\phi f(x,\xi )| \lesssim \eabs {x,\xi }^{N} \text ;
$$
\end{enumerate}
(Cf. \cite[Chapter 12]{Gc2}.)
\end{rem}

\par

\subsection{The Bargmann transform}\label{subsec1.3}
Next we recall some facts about the Bargmann transform. For every
$f\in \cS _{1/2} '(\rr d)$, the Bargmann transform $\mathfrak Vf$ is
the entire function on $\cc d$, defined by
\begin{equation}\label{bargdistrform}
(\mathfrak Vf)(z) =\scal f{\mathfrak A_d(z,\cdo )},
\end{equation}
where the Bargmann kernel $\mathfrak A_d$ is given by
$$
\mathfrak A_d(z,y)=\pi ^{-d/4} \exp \Big ( -\frac 12(\scal
zz+|y|^2)+2^{1/2}\scal zy\Big ).
$$
Here
$$
\scal zw = \sum _{j=1}^dz_jw_j,\quad 
z=(z_1,\dots ,z_d) \in \cc d,\quad  w=(w_1,\dots ,w_d)\in \cc d,
$$
and otherwise $\scal \cdo \cdo $ denotes the duality between test
function spaces and their corresponding duals.
We note that the right-hand side in \eqref{bargdistrform} makes sense
when $f\in \cS _{1/2} '(\rr d)$ and defines an element in the set $A(\cc d)$
of all entire functions on $\cc d$. In fact, $y\mapsto \mathfrak A_d(z,y)$
can be interpreted as an element
in $\cS _{1/2} (\rr d)$ with values in $A(\cc d)$.

\par

If in addition $f$ is an integrable function, then $\mathfrak Vf$ takes the form
$$
(\mathfrak Vf)(z) =\int \mathfrak A_d(z,y)f(y)\, dy,
$$
or
\begin{equation*}%\label{bargtransf}
(\mathfrak Vf)(z) =\pi ^{-d/4}\int _{\rr d}\exp \Big ( -\frac 12(\scal
z z+|y|^2)+2^{1/2}\scal zy \Big )f(y)\, dy .
\end{equation*}

\par

Several properties for the Bargmann transform were established
by Bargmann in \cite{B1,B2}. For example, in \cite{B1} it is proved
that $f\mapsto \mathfrak Vf$ is a bijective and
isometric map  from $L^2(\rr d)$ to the Hilbert space $A^2(\cc d)$, the set
of entire functions $F$ on $\cc  d$ which fulfill
\begin{equation}\label{A2norm}
\nm F{A^2}\equiv \Big ( \int _{\cc d}|F(z)|^2d\mu (z)  \Big )^{1/2}<\infty .
\end{equation}
Here $d\mu (z)=\pi ^{-d} e^{-|z|^2}\, d\lambda (z)$, where $d\lambda (z)$ is the
Lebesgue measure on $\cc d$, and the scalar product on $A^2(\cc d)$ is given by
\begin{equation}\label{A2scalar}
(F,G)_{A^2}\equiv  \int _{\cc d} F(z)\overline {G(z)}\, d\mu (z),\quad F,G\in
A^2(\cc d).
\end{equation}

\par

In \cite{B1} it is also proved that the Hermite functions are mapped by
the Bargmann transform into convenient monomials. More precisely,
for any multi-index $\alpha \in \nn d$, the Hermite function $h_\alpha$
of order $\alpha$ is defined by
$$
h_\alpha (x) = \pi ^{-d/4}(-1)^{|\alpha |}(2^{|\alpha |}\alpha
!)^{-1/2}e^{|x|^2/2}(\partial ^\alpha e^{-|x|^2}).
$$
It follows that
$$
h_{\alpha}(x)=  \frac{1}{(2\pi )^{d/2} \alpha !}
e^{-|x|^2/2}p_{\alpha}(x),
$$
for some polynomial $p_\alpha$ on $\rr d$, which is called the Hermite
polynomial of order $\alpha$.

\par

The set $\{ h_\alpha \} _{\alpha \in \nn d}$ is an orthonormal basis
for $L^2(\rr d)$. It is also a basis for any of the Gelfand-Shilov spaces
and their distribution spaces at above.

\par

In \cite{B1} it is then proved that
$$
(\mathfrak Vh_\alpha )(z) = \frac {z^\alpha}{\sqrt {\alpha !}},\qquad z\in \cc d .
$$

\par

Next we recall the links between the Bargmann transform and the short-time
Fourier transform, when the window function $\phi$ is given by
\begin{equation}\label{phidef}
\phi (x)=\pi ^{-d/4}e^{-|x|^2/2}.
\end{equation}
More precisely, let $S$ be the dilation operator given by
\begin{equation}\label{Sdef}
(SF)(x,\xi ) = F(2^{-1/2}x,-2^{-1/2}\xi ),
\end{equation}
when $F\in L^1_{loc}(\rr {2d})$. Then it
follows by straight-forward computations that
\begin{multline}\label{bargstft1}
  (\mathfrak{V} f)(z)  =  (\mathfrak{V}f)(x+i  \xi ) 
 =  (2\pi )^{d/2}e^{(|x|^2+|\xi|^2)/2}e^{-i\scal x\xi}V_\phi f(2^{1/2}x,-2^{1/2}\xi )
\\[1ex]
=(2\pi )^{d/2}e^{(|x|^2+|\xi|^2)/2}e^{-i\scal x\xi}(S^{-1}(V_\phi f))(x,\xi ),
\end{multline}
or equivalently,
\begin{multline}\label{bargstft2}
V_\phi f(x,\xi )  =  
(2\pi )^{-d/2} e^{-(|x|^2+|\xi |^2)/4}e^{-i \scal x \xi /2}(\mathfrak{V}f)
(2^{-1/2}x,-2^{-1/2}\xi).
\\[1ex]
=(2\pi )^{-d/2}e^{-i\scal x\xi /2}S(e^{-|\cdo |^2/2}(\mathfrak{V}f))(x,\xi ).
\end{multline}
For future references we observe that \eqref{bargstft1} and \eqref{bargstft2}
can be formulated into
\begin{equation}\label{BargStftlink}
\mathfrak V = U_{\mathfrak V}\circ V_\phi ,\quad \text{and}\quad
U_{\mathfrak V}^{-1} \circ \mathfrak V =  V_\phi ,
\end{equation}
where $U_{\mathfrak V}$ is the linear, continuous and bijective operator on
$\mathscr D'(\rr {2d})\simeq \mathscr D'(\cc d)$, given by
\begin{equation}\label{UVdef}
(U_{\mathfrak V}F)(x,\xi ) = (2\pi )^{d/2} e^{(|x|^2+|\xi |^2)/2}e^{-i\scal x\xi}
F(2^{1/2}x,-2^{1/2}\xi ) .
\end{equation}

\par

The next result shows that the image of the Bargmann transform of
the Gelfand-Shilov and tempered function and distribution spaces are
given by
\begin{equation*}
\begin{aligned}
\maclA _{\Sigma _t^s}(\cc d) &\equiv \sets {F\in A(\cc d)}{|F(z)|\lesssim
e^{|z|^2/2-\ep M_{s,t}(z)} \ \text{for every}\ \ep >0},
\\[1ex]
\maclA _{\maclS _t^s}(\cc d) &\equiv \sets {F\in A(\cc d)}{|F(z)|\lesssim
e^{|z|^2/2-\ep M_{s,t}(z)} \ \text{for some}\ \ep >0},
\\[1ex]
\maclA _{\mathscr S}(\cc d) &\equiv \sets {F\in A(\cc d)}{|F(z)|\lesssim
e^{|z|^2/2} \eabs z^{-N} \ \text{for every}\ N >0},
\\[1ex]
\maclA _{\mathscr S'}(\cc d) &\equiv \sets {F\in A(\cc d)}{|F(z)|\lesssim
e^{|z|^2/2} \eabs z^N
\ \text{for some}\ N >0},
\\[1ex]
\maclA _{(\maclS _t^s)'}(\cc d) &\equiv \sets {F\in A(\cc d)}{|F(z)|\lesssim
e^{|z|^2/2+\ep M_{s,t}(z)} \ \text{for every}\ \ep >0},
\\[1ex]
\maclA _{(\Sigma _t^s)'}(\cc d) &\equiv \sets {F\in A(\cc d)}{|F(z)|\lesssim
e^{|z|^2/2+\ep M_{s,t}(z)} \ \text{for some}\ \ep >0},
\end{aligned}
\end{equation*}
with canonical topologies. Here $M_{s,t}$ is given by
\begin{equation*}
M_{s,t}(x+i\xi ) = |x|^{1/t}+|\xi |^{1/s},\quad x,\xi \in \rr d.
\end{equation*}

\par

\begin{prop}\label{BargMapProp}
Let $s,t>1/2$, and let $\mathcal V_d$ be any of the spaces
\begin{alignat*}{3}
&\maclS _{1/2}(\rr d),&\quad &\maclS _t^{1/2}(\rr d),&\quad
&\maclS _{1/2}^s(\rr d), \quad \Sigma _t^s(\rr d),\quad \maclS _t^s(\rr d),
\quad \mathscr S(\rr d),
\\[1ex]
&\mathscr S'(\rr d),&\quad &(\maclS _t^s)'(\rr d),
&\quad &(\Sigma _t^s)'(\rr d).
\end{alignat*}
Then the map $f\mapsto \mathfrak V _d f$ is continuous and bijective
from $\mathcal V_d$ to $\maclA _{\mathcal V_d}(\cc d)$.
\end{prop}

\par

\begin{proof}
The result follows by a combination of Propositions \ref{stftGelfand2}
and \ref{stftGelfand2dist}, and Remark \ref{SchwFunctionSTFT},
with \eqref{BargStftlink}.
\end{proof}

\par

%%%%%%%%%%%%%%%%%%%%%%%
\section{Radial symmetric
elements and the Bargmann transform}\label{sec2}
%%%%%%%%%%%%%%%%%%%%%%%

\par

In this section we give a complete characterization for radial symmetric elements
under the Bargmann transform. The results are formulated for the broad
space $\cS '_{1/2}$ of Gelfand-Shilov  distributions. However, the applications
later on only involve elements in the smaller class $\sS '$. 

\par

We recall
that an element $f\in \cS '_{1/2}(\rr d)$ is called \emph{radial symmetric}, if the
pullback $U^*f$ is equal to $f$, for every unitary transformation $U$ on $\rr d$.
In the case when $f$ in addition is a measurable function,
then $f$ is radial symmetric, if and only if $f(x)=f_0(|x|)$ a.{\,}e., for some
measurable function $f_0$ on $[0,+\infty)$.

\par

\begin{thm}\label{Bargrotinv}
Let $f\in \cS '_{1/2}(\rr d)$. Then the following conditions are
equivalent:
\begin{enumerate}
\item $f$ is radial symmetric;

\vrum

\item if $U$ is unitary on $\rr d$, then
$(\mathfrak Vf)(Uz)=(\mathfrak Vf)(z)$;

\vrum

\item $(\mathfrak Vf)(z) =F_0(\scal zz )$, for some entire function $F_0$
on $\mathbf C$.

\vrum

\item $(f,h_\alpha )=0$ for every $\alpha =(\alpha _1,\dots ,\alpha _d)\in \nn d$
such that at least one of $\alpha _j$ is odd, and
$$
\frac {\alpha !}{\sqrt {(2\alpha )!}} (f,h_{2\alpha})
= \frac {\beta !}{\sqrt {(2\beta )!}} (f,h_{2\beta})
$$
%$$
%{\alpha !}((2\alpha )!)^{-1/2} (f,h_{2\alpha})
%= {\beta !}((2\beta )!)^{-1/2} (f,h_{2\beta})
%$$
for any $\beta \in \nn d$ with $|\alpha| = |\beta|$.

\end{enumerate}
\end{thm}

\par

Here $Uz$ is defined as $Ux+iUy$, when $z=x+iy$ and $x,y\in \rr d$.

\par

\begin{proof}
Let $U$ be unitary on $\rr d$ and let
$f\in \cS '_{1/2}(\rr d)$. The equivalence between (1) and (2) follows
if we prove that
\begin{equation}\label{unitpullbackBargm}
(\mathfrak V(U^*f))(z) = (\mathfrak Vf)(Uz).
\end{equation}

\par

Let $U^t$ be the transpose of $U$. Then
\begin{multline*}
(\mathfrak V(U^*f))(z) = \pi ^{-d/4} \Big \langle f,\exp \Big ( -\frac 12(\scal
zz+|U^t \cdo |^2)+2^{1/2}\scal z{U^t\cdo }\Big )\Big \rangle
\\[1ex]
= \pi ^{-d/4} \Big \langle f,\exp \Big ( -\frac 12(\scal
{Uz}{Uz}+|\cdo |^2)+2^{1/2}\scal {Uz}{\cdo }\Big )\Big \rangle
=(\mathfrak Vf)(Uz).
\end{multline*}
Here the second equality follows from the fact that $\scal {Uz}{Uw}=\scal zw$,
when $z,w\in \cc d$. This gives \eqref{unitpullbackBargm}.

\par

Next assume that (1) holds, and let $z=x\in \rr d \subseteq \cc d$
be fixed. Let $U$ be unitary such that $Ux=|x|e_1$, where $e_1$ is the first unit
vector in $\rr d$.

\par

Since $U^tf=f$, $U^t=U^{-1}$ and $|\det U|=1$, we get
\begin{multline*}
(\mathfrak Vf)(x) = \pi ^{-d/4}\Big \langle f,\exp \Big ( -\frac 12(|x|^2+|\cdo |^2)
+2^{1/2}\scal x\cdo \Big )\Big \rangle
\\[1ex]
= \pi ^{-d/4}\Big \langle U^tf,\exp \Big ( -\frac 12(|x|^2+|\cdo |^2)
+2^{1/2}\scal {x}\cdo \Big )\Big \rangle
\\[1ex]
= \pi ^{-d/4}\Big \langle f,\exp \Big ( -\frac 12(|x|^2+|U\cdo |^2)
+2^{1/2}\scal {Ux}\cdo \Big )\Big \rangle .
\end{multline*}

\par

Now we have
$$
\exp \Big ( -\frac 12(|x|^2+|Uy |^2) +2^{1/2}\scal {Ux}y \Big )
=
\exp \Big ( -\frac 12(|x|^2+|y|^2) +2^{1/2}|x|y_1 \Big ).
$$
Hence if
$$
E(t,y) = \pi ^{-d/4}\exp \Big ( -\frac 12(t^2+|y|^2) +2^{1/2}ty_1 \Big ),
\quad t\in \mathbf R,\ y\in \rr d,
$$
then
\begin{equation}\label{BargmannEfunct}
(\mathfrak Vf)(x) = \scal f{E(|x|,\cdo )}.
\end{equation}

\par

Since similar facts hold if we should have chosen $U$ such that $Ux=-|x|e_1$,
it also follows that
\begin{equation}\tag*{(\ref{BargmannEfunct})$'$}
(\mathfrak Vf)(x) = \scal f{E(-|x|,\cdo )}.
\end{equation}
By taking the mean-value of \eqref{BargmannEfunct} and
\eqref{BargmannEfunct}$'$, we get
$$
(\mathfrak Vf)(x) = \scal f{E_0(|x|^2,\cdo )},
$$
when
$$
E_0(t^2,y) = \pi ^{-d/4}\exp \Big ( -\frac 12(t^2+|y|^2)\Big ) \cosh (2^{1/2}ty_1),
$$
or equivalently, when
$$
E_0(s,y) = \pi ^{-d/4}\exp \Big ( -\frac 12(s+|y|^2)\Big ) \sum _{k=0}^\infty
\frac {(2y_1^2s)^k}{(2k)!},\quad s\ge 0.
$$

\par

Evidently, $E_0$ is uniquely extendable to an entire function on $\mathbf C
\times \cc d$, and if
$$
F_0(w)\equiv \scal f{E_0(w,\cdo )},\quad w\in \mathbf C,
$$
then $F_0$ is entire and $(\mathfrak Vf)(x) = F_0(|x|^2)$.
Now, since $(\mathfrak Vf)(z)$ and $F_0(\scal zz)$ are entire functions
on $\cc d$, which coincide on $\rr d$, it follows that $(\mathfrak Vf)(z)=
F_0(\scal zz)$ for all $z\in \cc d$, and (3) follows.

\par

Next we prove that (3) implies (1). Therefore assume that (3) holds, let
$U$ be unitary on $\rr d$, and recall that
$\scal {Uz}{Uz}=\scal zz$. Hence \eqref{unitpullbackBargm} gives
$$
(\mathfrak V(U^*f))(z) = (\mathfrak Vf)(Uz) = F_0(\scal {Uz}{Uz})
= F_0(\scal zz) = (\mathfrak Vf)(z).
$$
Since the Bargmann transform is injective on $\cS '_{1/2}$, it follows
that $U^*f=f$, and (1) holds.

\par

To conclude the proof it is sufficient to prove that (3) is equivalent to (4).

\par

Let $f\in \mathcal V_d$, $F=\mathfrak Vf$, and let
$$
H_\alpha (z)\equiv \frac {z^\alpha}{\sqrt {\alpha !}}.
$$
Then $\mathfrak Vh_\alpha =H_\alpha$, and $F$ possess
the unique expansion
$$
F = \sum _{\alpha \in \nn d} a_\alpha H_\alpha,
$$
where
$$
a_\alpha = (F,H_\alpha )_{A^2} = (f,h_\alpha )_{L^2}.
$$

\par

The assertion (4) is true, if and only if $F$ can be written as
$$
F(z) = \sum _{k=0}^\infty c_k \left (\sum
_{\gamma \in \nn d,\ |\gamma |=k} \frac {z^{2\gamma}}{\gamma !}\right ).
$$
It now follows from the uniqueness of the expansions and the identities above
that the latter identity holds, if and only if $a_\alpha =0$ when at least
one $\alpha _j$ is odd, and that
$$
\frac {\alpha !}{\sqrt {(2\alpha )!}} (f,h_{2\alpha})
$$
only depends on $|\alpha |$. The proof is complete.
\end{proof}

\par

We shall use the previous result to show that there is a way to
relate any radial symmetric element in $\mathcal V_d$ in Proposition
\ref{BargMapProp} with an element in $\mathcal V_1$. More precisely, let
$\mathcal V_d$ be as in Proposition \ref{BargMapProp}, $f\in
\mathcal V_d$ be radial symmetric, and let $\mathfrak V_df(z)=F_0(\scal zz)$,
where $F_0$ is the same as in Theorem \ref{Bargrotinv}. Then
$$
{\mathbf C} \owns z\mapsto F_0(z^2)
$$
is entire, and belongs to $\mathcal V_1$, in view of Proposition
\ref{BargMapProp}. By using Proposition \ref{BargMapProp} again
we get the following result.

\par

\begin{prop}\label{ReductionProp}
Let $\mathcal V_d$ be as in Proposition \ref{BargMapProp}, $f\in
\mathcal V_d$ be radial symmetric, and let $\mathfrak V_df(z)=F_0(\scal zz)$,
where $F_0$ is the same as in Theorem \ref{Bargrotinv}. Then
there is a unique $f_0\in \mathcal V_1$ such that
$$
\mathfrak V_1f_0(z)=F_0(z^2),\quad z\in \mathbf C.
$$
\end{prop}

\par

\begin{rem}
Let $f$ and $F_0$ be the same as in Proposition \ref{ReductionProp}.
Then it follows from Proposition \ref{BargMapProp} that $F_0$
belongs to the subclass
$$
\mathcal A_0(\mathbf C) \equiv
\sets {F\in A(\mathbf C)}{|F(z)|\lesssim
e^{\ep |z|} \ \text{for some}\ \ep >0}
$$
of $\maclA _{\maclS _{1/2}}(\mathbf C)$. Hence, $F_0 = \mathfrak V_1f_1(z)$
for some $f_1\in \maclS _{1/2}(\mathbf R)$.

\par

It seems to be difficult to assign to the functions or distributions
$f_1$ and $f_0$ any specific roles. For example, if $f$ in addition
is a function, then $f(x) \neq cf_j(|x|)$, $j=0,1$, for every constant
$c\in \mathbf C$. This follows by straight-forward control of the
Taylor expansions of $F_0(z^2)$, $\mathfrak V_1f_1(z)$
and $\mathfrak V_1f_0(z)$.
\end{rem}

\end{document}